\documentclass[graybox]{svmult}

\usepackage{helvet}         
\usepackage{courier}        
\usepackage{type1cm}        
\usepackage{makeidx}         
\usepackage{graphicx}        
\usepackage{multicol}        
\usepackage[bottom]{footmisc}
\usepackage{dynkin-diagrams}
\usepackage{tikz-cd}
\newcommand{\ve}{\varepsilon}

\usepackage{amsmath}
\usepackage{amsxtra}
\usepackage{amscd}
\usepackage{amsfonts}
\usepackage{amssymb}
\usepackage{color}
\usepackage{ulem}
\begin{document}

\title*{Non-local integrals of motion for deformed
$W$-algebras of types $A_l, D_l, E_{6,7,8}$}

\titlerunning{Non-local Integrals of Motion of types $A_l, D_l, E_{6,7,8}$} 
\author{Michio Jimbo and Takeo Kojima}
\institute{Michio Jimbo\at 
Emeritus professor, Rikkyo University, Toshima-ku, Tokyo 171-8501, Japan, \email{
jimbomm@rikkyo.ac.jp}
\and Takeo Kojima\at 
Yamagata University, Jonan 4-3-16, Yonezawa 992-8510, Japan
\email{
kojima@yz.yamagata-u.ac.jp}}

\maketitle

\noindent
\textit{To the memory of Masatoshi Noumi}
\vspace{10mm}

\noindent
\bigskip

\abstract{
We present
an infinite set of non-local integrals of motion for 
deformed $W$-algebras of types $A_l, D_l$, and $E_{6,7,8}$.
They can be regarded as a two-parameter deformation of trace of the monodromy matrix of the ${\mathfrak g}$-KdV theory. 
Commutativity of the non-local integrals of motion is shown in the case of
$A_l$ and $D_l$ by a direct calculation. In the case of
$E_{6,7,8}$ it is a conjecture. 
}

\section{Introduction}
\label{sec:1}

First let us recall the construction of integrals of motion
in soliton theory, taking the KdV equation as an example.
The KdV equation can be written in the Poisson bracket form 
\cite{Gardner, Faddeev-Zakharov}
\begin{eqnarray}
u_t=\{u, H\}=6u u_x+u_{x x x},~~~H=\int u(x)^2 dx.
\nonumber
\end{eqnarray}
Here the Poisson bracket is given by
\begin{eqnarray}
\{u(x), u(y)\}=
\left(\frac{1}{2}\partial_x^3+2u(x)\partial_x+u_x(x)\right)\delta(x-y).
\nonumber
\end{eqnarray}
This defines a Poisson algebra 
called the classical $W$-algebra ${\mathcal W}(A_1)$
(aka classical Virasoro algebra).

The KdV equation has the Lax representation \cite{Lax1}
in terms of  $L=\partial_x^2+u$ and 
$M=4 \partial_x^3+6u \partial_x+3 u_x$:
\begin{eqnarray}
\partial_{t} L=[M, L]=6u u_x+u_{x x x}.\nonumber
\end{eqnarray}
Let $u(x,t)$
be a solution to the KdV equation satisfying the
periodic boundary condition $u(x+2\pi,t)=u(x,t)$.
Let $\psi_1(x,t), \psi_2(x,t)$ be two linearly independent solutions of the differential equation:
$L \psi_j(x,t)=\lambda \psi_j(x,t)$.
Introduce the transfer matrix
\begin{eqnarray}
{\bf T}^{(cl)}(\lambda)={\rm tr} \left({\mathbf M}^{(cl)}(\lambda)\right),
\nonumber
\end{eqnarray}
where ${\bf M}^{(cl)}(\lambda)$ is
the monodromy matrix of solutions
$\psi_1(x,t), \psi_2(x,t)$:
\begin{eqnarray}
(\psi_1(x+2\pi,t), \psi_2(x+2\pi,t))=(\psi_1(x,t), \psi_2(x,t)){\mathbf M}^{(cl)}(\lambda).
\nonumber
\end{eqnarray}
The classical inverse scattering method
\cite{Lax2, Date-Tanaka, Marchenko, Faddeev-Takhtajan}
yields
\begin{eqnarray}
\{{\bf T}^{(cl)}(\lambda), {\bf T}^{(cl)}(\mu)\}=0,~~~
\partial_{t} {\bf T}^{(cl)}(\lambda)=0.
\nonumber
\end{eqnarray}
Hence the Taylor coefficients of
${\displaystyle {\mathbf T}^{(cl)}(\lambda)
=\sum_{n=0}^\infty {\mathbf G}_n^{(cl)} \lambda^{n}}$
are independent of $t$
and are mutually in involution,
\begin{eqnarray}
\{{\mathbf G}_m^{(cl)}, {\mathbf G}_n^{(cl)}\}=0,~~~
\partial_{t} {\mathbf G}_{n}^{(cl)}=0,~~~
m,n=1,2,3,\ldots.
\nonumber
\end{eqnarray}
We call $\{{\bf G}_n^{(cl)}\}_{n=1}^\infty$
classical non-local integrals of motion.

The construction above was generalized by
Drinfeld-Sokolov \cite{Drinfeld-Sokolov} for 
an arbitrary affine Lie algebra. \footnote{
While $W$ algebras are associated with simple Lie algebras,
integrable hierarchies correspond to affine Lie algebras.
}
They introduced the $\mathfrak{g}$-KdV theory, 
the classical $W$-algebra ${\mathcal W}({\mathfrak g})$
and the classical non-local integrals of motion associated with it.

The $W$-algebra ${\mathcal W}_k({\mathfrak g})$ is a one-parameter deformation of the classical $W$-algebra
${\mathcal W}({\mathfrak g})$ and
plays a role of symmetry underpinning
conformal field theory (CFT)
\cite{Feigin-Frenkel1, Kac-Roan-Wakimoto, Feigin-Frenkel2, Belavin-Polyakov-Zamolodchikov}.
Bazhanov et al. applied the quantum inverse scattering method to
the CFT analogue of the KdV theory and
obtained remarkable results
\cite{
Bazhanov-Lukyanov-Zamolodchikov1,
Bazhanov-Lukyanov-Zamolodchikov2,
Bazhanov-Lukyanov-Zamolodchikov3,
Bazhanov-Lukyanov-Zamolodchikov4, Dorey-Tateo,
Bazhanov-Hibberd-Khoroshkin,Feigin-Frenkel3}.
The deformed $W$-algebra 
${\mathcal W}_{q,t}({\mathfrak g})$
is a two-parameter deformation of the classical $W$-algebra ${\mathcal W}({\mathfrak g})$
\cite{Shiraishi-Kubo-Awata-Odake,
Awata-Kubo-Odake-Shiraishi,
Feigin-Frenkel4,
Sevostyanov, 
Frenkel-Reshetikhin,
Feigin-Jimbo-Mukhin-Vilkoviskiy}. 
The deformed $W$ algebra ${\mathcal W}_{q,t}({\mathfrak g})$
is an interesting object for its own sake 
as well as its connection with other topics, including
quantum toroidal algebras, Macdonald polynomials, and their generalizations.

In this paper we present
an infinite set of non-local integrals of motion for the deformed $W$-algebra of types $A_l, D_l$, and $E_{6,7,8}$.
They can be regarded as a two-parameter deformation of the non-local integrals of motion in the ${\mathfrak g}$-KdV equation. 
Commutativity of the non-local integrals of motion is shown in the case of $A_l$ and $D_l$ by a
direct calculation
\cite{Feigin-Kojima-Shiraishi-Watanabe, Kojima-Shiraishi, Feigin-Jimbo-Mukhin1, Jimbo-Kojima, Feigin-Jimbo-Mukhin2}.
In the case of $E_{6,7,8}$ it is a conjecture.

The plan of the paper is as follows.
In Section \ref{sec:2}
we set up our notation.
In Section \ref{sec:3}
we introduce the screening currents
for the deformed $W$-algebra of types $A_l, D_l$ and $E_{6,7,8}$.
In Section \ref{sec:4}  we give free field constructions of the non-local integrals of motion
for the deformed $W$-algebra.
We state a theorem and a conjecture.
In Section \ref{sec:5} we give discussions.

\section{Preliminaries}

\label{sec:2}

In this section we introduce our notation.

\subsection{Notation}

Throughout this paper, we fix a real number $\beta$ with $\beta>2$, and a complex number $q$ with $0<|q|<1$.
We fix also
a complex number $\tau$ with ${\rm Im}(\tau)>0$ such that
$q=e^{-2\pi \sqrt{-1}/\tau}$.

We use symbols for infinite products,
\begin{eqnarray}
(a;q)_\infty=\prod_{k=0}^\infty (1-a q^k),~~~
(a_1,a_2,\ldots,a_N;q)_\infty=\prod_{i=1}^N (a_i;q)_\infty
\nonumber
\end{eqnarray}
for complex numbers $a, a_1, \ldots, a_N$.
We set
\begin{eqnarray}
\Theta_q(z)=(q,w,qw^{-1};q)_\infty.\nonumber
\end{eqnarray}
The symbol $\theta(u)$ stands for the elliptic theta function
satisfying
\begin{eqnarray}
\theta(u+1)=-\theta(u)=\theta(-u),~~~
\theta(u+\tau)=-e^{-2\pi\sqrt{-1} u-\pi \sqrt{-1}\tau}\theta(u).\nonumber
\end{eqnarray}
Explicitly it is given by
\begin{eqnarray}
\theta(u)=q^{\frac{u^2}{2}-\frac{u}{2}}\Theta_q(q^u).\nonumber
\end{eqnarray}

We use the symbols for
lexicographic ordered product
\begin{eqnarray}
&
\overset{\longrightarrow}{\displaystyle
\prod_{1\leq i \leq N}}A_i=A_1 A_2 \cdots A_N,
\nonumber
\\
&
\overset{\longrightarrow}{\displaystyle
\prod_{1\leq i \leq M}}
\overset{\longrightarrow}{\displaystyle
\prod_{1\leq j \leq N}}A_{i,j}=
\overset{\longrightarrow}{
\displaystyle
\prod_{1\leq j \leq N}}A_{1,j} 
\overset{\longrightarrow}{
\displaystyle
\prod_{1\leq j \leq N}}A_{2,j} 
\cdots
\overset{\longrightarrow}{
\displaystyle
\prod_{1\leq j \leq N}}A_{M,j}.\nonumber
\end{eqnarray}

\subsection{Simply-laced Lie algebras}

Let $\mathfrak{g}$ be a simply laced complex simple
Lie algebra of rank $l$, 
with simple roots $\alpha_i$, $1\le i\le l$. 
Let further $\alpha_0$ be the negative of the highest root.

Denote by $A^{[{\mathfrak g}]}=\bigl(A_{i,j}^{[{\mathfrak g}]}\bigr)_{i,j=0}^l$
the extended Cartan matrix with entries
$A_{i,j}^{[{\mathfrak g}]}=(\alpha_i|\alpha_j)$, where $(~|~)$ is the 
standard bilinear form on $\mathfrak{g}$
such that $(\alpha_i|\alpha_i)=2$. 
The Kac numbers
$a_i^{[{\mathfrak g}]}$, $0\le i\le l$, 
are relatively prime positive integers
satisfying  $\sum_{i=0}^l a_i^{[{\mathfrak g}]}\alpha_i=0$. 

We give below these data for each $\mathfrak{g}$. 
In the sequel 
 ${\bf R}^{n}$ stands for 
an $n$-dimensional real Euclidean space equipped
with the standard basis $\ve_1, \ve_2,\ldots, \ve_n$ and the bilinear form:
$(\ve_i|\ve_j)=\delta_{i,j}$.

\subsubsection{$A_l$, $l\geq 1$}

We have $n=l+1\geq 2$,
\begin{eqnarray}
&\alpha_0=-\ve_1+\ve_{l+1},
\alpha_1=\ve_1-\ve_2,\ldots, 
\alpha_i=\ve_i-\ve_{i+1},\ldots,
\alpha_l=\ve_l-\ve_{l+1},
\nonumber
\end{eqnarray}
and the Kac numbers are 
\begin{eqnarray}
&(a_0^{[{A_l}]},a_1^{[{A_l}]},\ldots, a_l^{[{A_l}]})=(\overbrace{1,1,1,\ldots,1}^{l+1}).
\nonumber
\end{eqnarray}

\subsubsection{$D_l$, $l\geq 4$}

We have $n=l\geq 4$,
\begin{eqnarray}
&&\alpha_0=-\ve_1-\ve_{2},
\alpha_1=\ve_1-\ve_2, \ldots, \alpha_i=\ve_i-\ve_{i+1},
\nonumber\\
&&
\ldots, \alpha_{l-1}=\ve_{l-1}-\ve_l,
\alpha_l=\ve_{l-1}+\ve_{l} {\color{blue} ,}
\nonumber
\end{eqnarray}
and the Kac numbers are 
\begin{eqnarray}
&
(a_0^{[{D_l}]},a_1^{[{D_l}]},\ldots,a_l^{[{D_l}]})
=(1,1,\overbrace{2,2,\ldots,2,2}^{l-3},1,1)\,.
\nonumber
\end{eqnarray}

\subsubsection{
$E_l$, $l=6,7,8$}
~\\
$\bullet$~$E_6$: 
We have $n=7$,
\begin{eqnarray}
&&
\alpha_0=\frac{1}{2}(\ve_6+\ve_7)-\frac{1}{2}(\ve_1+\ve_2+\ve_3+\ve_4+\ve_5+\ve_8),
\nonumber\\
&&
\alpha_1=\frac{1}{2}(\ve_1+\ve_8)-\frac{1}{2}(\ve_2+\ve_3+\ve_4+\ve_5+\ve_6+\ve_7),
\nonumber\\
&&
\alpha_2=\ve_1+\ve_2,
\alpha_3=\ve_2-\ve_1,
\alpha_4=\ve_3-\ve_2,\nonumber
\\
&&
\alpha_5=\ve_4-\ve_3,
\alpha_6=\ve_5-\ve_4\,,
\nonumber
\end{eqnarray}
and the Kac numbers are
\begin{eqnarray}
&(a_0^{[{E_6}]},a_1^{[{E_6}]},\ldots, a_6^{[{E_6}]})
=(1,1,2,2,3,2,1)\,.
\nonumber
\end{eqnarray}
~\\
$\bullet$~$E_7$: 
We have $n=8$,  
\begin{eqnarray}
&&
\alpha_0=\ve_7-\ve_8,\nonumber\\
&&
\alpha_1=
\frac{1}{2}(\ve_1+\ve_8)-\frac{1}{2}(\ve_2+\ve_3+\ve_4+\ve_5+\ve_6+\ve_7),
\nonumber\\
&&
\alpha_2=\ve_1+\ve_2,
\alpha_3=\ve_2-\ve_1,
\alpha_4=\ve_3-\ve_2,\nonumber
\\
&&
\alpha_5=\ve_4-\ve_3,
\alpha_6=\ve_5-\ve_4,
\alpha_7=\ve_6-\ve_5\,,
\nonumber
\end{eqnarray}
and the Kac numbers are
\begin{eqnarray}
&(a_0^{[{E_7}]},a_1^{[{E_7}]},\ldots,a_7^{[{E_7}]})
=(1,2,2,3,4,3,2,1)\,.
\nonumber
\end{eqnarray}
~\\
$\bullet$~$E_8$: 
We have $n=8$,
\begin{eqnarray}
&&\alpha_0=-\ve_7-\ve_8,\nonumber\\
&&\alpha_1=\frac{1}{2}(\ve_1+\ve_8)
-\frac{1}{2}(\ve_2+\ve_3+\ve_4+\ve_5+\ve_6+\ve_7),\nonumber\\
&&
\alpha_2=\ve_1+\ve_2,
\alpha_3=\ve_2-\ve_1,
\alpha_4=\ve_3-\ve_2,
\alpha_5=\ve_4-\ve_3,
\nonumber\\
&&
\alpha_6=\ve_5+\ve_4,
\alpha_7=\ve_6-\ve_5,
\alpha_8=\ve_7-\ve_6\,,
\nonumber
\end{eqnarray}
and the Kac numbers are
\begin{eqnarray}
&(a_0^{[{E_8}]},a_1^{[{E_8}]},\ldots, a_8^{[{E_8}]})
=(1,2,3,4,6,5,4,3,2)\,.
\nonumber
\end{eqnarray}

\section{Screening currents}

\label{sec:3}

In this section we recall 
the screening currents $S_i^{[{\mathfrak g}]}(z)$, $0\le i\le l$,
associated with
${\mathfrak g}=A_l, D_l$ and $E_{6,7,8}$
\cite{Frenkel-Reshetikhin, 
Feigin-Jimbo-Mukhin-Vilkoviskiy,
Feigin-Kojima-Shiraishi-Watanabe, 
Kojima-Shiraishi}.

\subsection{Heisenberg algebra}

Define the extended
$q$-Cartan matrix $A^{[{\mathfrak g}]}(q,\beta)=(A_{i,j}^{[{\mathfrak g}]}(q,\beta))_{i,j=0}^l$
of types ${\mathfrak g}=A_l$, $l \geq 1$, $D_l$, $l \geq 4$, 
and $E_{l}$, $l=6,7,8$, 
as follows.
\begin{eqnarray}
A_{i,j}^{[{\mathfrak g}]}(q,\beta)=
\frac{q^{\frac{\beta-1}{2}}-q^{\frac{1-\beta}{2}}}{q^{\frac{1}{2}}-q^{-\frac{1}{2}}}
\times
\left\{\begin{array}{cc}
q^{\frac{\beta}{2}}+q^{-\frac{\beta}{2}},
& A_{i,j}^{[{\mathfrak g}]}=2,\\
-1,
&A_{i,j}^{[{\mathfrak g}]}=-1,\\
(-1) (
q^{\frac{\beta}{2}-\gamma}+q^{\gamma-\frac{\beta}{2}}),
&A_{i,j}^{[{\mathfrak g}]}=-2,\\
0,&A_{i,j}^{[{\mathfrak g}]}=0,
\end{array}
\right.~~~0\leq i,j \leq l.
\nonumber
\end{eqnarray}
In the case of $A_{i,j}^{[A_1]}=-2$,
we add an additional parameter $\gamma$ with $0<\gamma<\beta$ to remove
``a singularity at $\gamma=\beta$" 
in the free field construction of the non-local integrals of motion \cite{Feigin-Kojima-Shiraishi-Watanabe}.

Let ${\mathcal H}^{[{\mathfrak g}]}$ be 
the Heisenberg algebra
generated by
elements $B_{i,m}^{[{\mathfrak g}]}$, $0\leq i \leq l, m \in {\mathbf Z}_{\neq 0}$, and
the zero mode operators $P_{\alpha_i}^{[{\mathfrak g}]}$, $Q_{\alpha_i}^{[{\mathfrak g}]}$, $0\leq i \leq l$,
with relations
\begin{eqnarray}
&&
[B_{i,m}^{[{\mathfrak g}]}, B_{j,n}^{[{\mathfrak g}]}]=\frac{1}{m} A_{i,j}^{[{\mathfrak g}]}(q^m, \beta)\delta_{m+n,0},~~~0\leq i,j \leq l, m,n \in {\bf Z}_{\neq 0},
\nonumber\\
&&
[P_{\alpha_i}^{[{\mathfrak g}]},Q_{\alpha_j}^{[{\mathfrak g}]}]=(1-\beta)A_{i,j}^{[{\mathfrak g}]},~~~0 \leq i,j \leq l.\nonumber
\end{eqnarray}
We use the normal ordering $:~:$ that satisfies
\begin{eqnarray}
:B_{i,m}^{[\mathfrak{g}]} B_{j,n}^{[{\mathfrak g}]}:=\left\{
\begin{array}{cc}
B_{i,m}^{[\mathfrak{g}]} B_{j,n}^{[{\mathfrak g}]},& m<0,\\
B_{j,n}^{[\mathfrak{g}]} B_{i,m}^{[{\mathfrak g}]},& m>0,
\end{array}
\right.~~~
:P_{\alpha_i}^{[{\mathfrak g}]} Q_{\alpha_j}^{[{\mathfrak g}]}:=:Q_{\alpha_j}^{[{\mathfrak g}]}P_{\alpha_i}^{[{\mathfrak g}]}:=
Q_{\alpha_j}^{[{\mathfrak g}]}P_{\alpha_i}^{[{\mathfrak g}]},~~~0\leq i,j \leq l.
\nonumber
\end{eqnarray}
The zero mode operators $P_{\alpha_i}^{[{\mathfrak g}]}$, $Q_{\alpha_i}^{[{\mathfrak g}]}$, $0\leq i \leq l$ are linearly dependent.
\begin{eqnarray}
\sum_{i=0}^l a_i^{[{\mathfrak g}]}P_{\alpha_i}^{[{\mathfrak g}]}=0,\quad
\sum_{i=0}^l a_i^{[{\mathfrak g}]}Q_{\alpha_i}^{[{\mathfrak g}]}=0.
\notag
\end{eqnarray}
Fix
parameters 
$\mu_0^{[{\mathfrak g}]},\ldots, \mu_l^{[{\mathfrak g}]} \in {\mathbf C}$ such that
\begin{eqnarray}
\sum_{i=0}^l a_i^{[{\mathfrak g}]} \mu_i^{[{\mathfrak g}]}=0.
\notag
\end{eqnarray}
Let $|\mu^{[{\mathfrak g}]}\rangle\neq 0$ be the Fock vacuum of the Fock space $\pi_{\mu^{[{\mathfrak g}]}}$ of the Heisenberg algebra
${\mathcal H}^{[{\mathfrak g}]}$
such that 
\begin{align}
B_{i,m}^{[{\mathfrak g}]} |\mu^{[{\mathfrak g}]}\rangle=0,\quad
P_{\alpha_i}^{[{\mathfrak g}]}|\mu^{[{\mathfrak g}]} \rangle
=\mu_i^{[{\mathfrak g}]} |\mu^{[{\mathfrak g}]}\rangle,\quad 0\leq i \leq l, m>0.\notag
\end{align}
Later, we give non-local integrals of motion acting on the Fock space $\pi_{\mu^{[{\mathfrak g}]}}$.

\subsection{Screening currents}

Define the screening currents $S_i^{[{\mathfrak g}]}(z)$ for the deformed $W$-algebra associated with
${\mathfrak g}=A_l$, $D_l$, and $E_{6,7,8}$
as follows.
\begin{eqnarray}
S_i^{[{\mathfrak g}]}(z)=
z^{1-\beta} e^{Q_{\alpha_i}^{[{\mathfrak g}]}}z^{P_{\alpha_i}^{[{\mathfrak g}]}}
:\exp\left(\sum_{n \neq 0}B_{i,n}^{[{\mathfrak g}]}z^{-n}\right):,~~~0\leq i \leq l.\nonumber
\end{eqnarray}
The screening currents $S_i^{[{\mathfrak g}]}(z)$, $0\leq i \leq l$, satisfy
\begin{eqnarray}
\theta(u-v+\beta)\theta(u-v)S_i^{[{\mathfrak g}]}(z)S_j^{[{\mathfrak g}]}(w)
&=&
\theta(v-u+\beta)\theta(v-u)S_j^{[{\mathfrak g}]}(w)S_i^{[{\mathfrak g}]}(z),~~~A_{i,j}^{[{\mathfrak g}]}=2,
\nonumber
\\
\theta\left(u-v+\frac{\beta}{2}\right)^{-1}S_i^{[{\mathfrak g}]}(z)
S_j^{[{\mathfrak g}]}(w)
&=&
\theta\left(v-u+\frac{\beta}{2}\right)^{-1}
S_j^{[{\mathfrak g}]}(w)S_i^{[{\mathfrak g}]}(z),~~~A_{i,j}^{[{\mathfrak g}]}=-1,\nonumber\\
S_i^{[{\mathfrak g}]}(z)
S_j^{[{\mathfrak g}]}(w)&=&
S_j^{[{\mathfrak g}]}(w)S_i^{[{\mathfrak g}]}(z),~~~A_{i,j}^{[{\mathfrak g}]}=0.
\nonumber
\end{eqnarray}

\section{Non-local integrals of motion}
\label{sec:4}

In this section we give a free field construction
of the non-local integrals of motion for the deformed $W$-algebra 
of types $A_l, D_l$, and $E_{6,7,8}$.

\subsection{Free field construction}

We fix parameters 
$\mu_0^{[{\mathfrak g}]},\ldots, \mu_l^{[{\mathfrak g}]}$ such that
$\sum_{i=0}^l a_i^{[{\mathfrak g}]} \mu_i^{[{\mathfrak g}]}=0$, which define the action of the Heisenberg algebra ${\mathcal H}^{[{\mathfrak g}]}$
on the Fock space $\pi_{\mu^{[{\mathfrak g}]}}$.
Let $X^{[\mathfrak{g}]}$ be the space of entire functions
$\vartheta^{[{\mathfrak g}]}(v_0,\ldots, v_l)$ 
satisfying
\begin{eqnarray}
&&\vartheta^{[{\mathfrak g}]}(v_0,\ldots, v_i+1,\ldots, v_l)=\vartheta^{[{\mathfrak g}]}(v_0,\ldots, v_i,\ldots, v_l),~~~0\leq i \leq l,\nonumber\\
&&
\vartheta^{[{\mathfrak g}]}(v_0,\ldots, v_i+\tau,\ldots, v_l)=
e^{-2\pi \sqrt{-1}(\sum_{j=0}^l A_{i,j}^{[{\mathfrak g}]}v_j-\mu_i^{[{\mathfrak g}]}
+\tau)}\vartheta^{[{\mathfrak g}]}(v_0,\ldots, v_i,\ldots, v_l),~~~0\leq i \leq l.\nonumber
\end{eqnarray}

Consider first the case $\mathfrak{g}\neq A_1$.
For each $N\ge0$ and $\vartheta^{[{\mathfrak g}]}\in 
X^{[\mathfrak{g}]}$ 
we define
${\bf G}_N^{[{\mathfrak g}]}(\vartheta^{[{\mathfrak g}]})$ as follows:
\begin{eqnarray}
{\bf G}_N^{[{\mathfrak g}]}(\vartheta^{[{\mathfrak g}]})
&=&
\prod_{i=0}^l \prod_{a=1}^{N a_i^{[{\mathfrak g}]
}} \int_{|z_{i,a}|=1}\frac{dz_{i,a}}{2\pi \sqrt{-1}z_{i,a}}
\overset{\longrightarrow}{\displaystyle \prod_{0\leq i \leq l}} 
\overset{\longrightarrow}{\prod_{~1\leq a \leq N a_i^{[{\mathfrak g}]}}}S_i^{[{\mathfrak g}]}(z_{i,a})
\label{def:non-local1}
\\
&\times&
\frac{\displaystyle \prod_{i=0}^l \prod_{1\leq a<b \leq N a_i^{[{\mathfrak g}]}}
\theta(u_{i,a}-u_{i,b}+\beta)
\theta(u_{i,a}-u_{i,b})}{\displaystyle
\prod_{0\leq i<j \leq l
\atop{A_{i,j}^{[{\mathfrak g}]}=-1}}
\prod_{a=1}^{N a_i^{[{\mathfrak g}]}}
\prod_{b=1}^{N a_j^{[{\mathfrak g}]}}
\theta\left(u_{i,a}-u_{j,b}+\frac{\beta}{2}\right)}
\times
\vartheta^{[{\mathfrak g}]}
\left(\sum_{a=1}^{N a_0^{[{\mathfrak g}]}}u_{0,a},\ldots,
\sum_{a=1}^{N a_i^{[{\mathfrak g}]}}u_{i,a},\ldots,
\sum_{a=1}^{N a_l^{[{\mathfrak g}]}}u_{l,a}
\right).\nonumber
\end{eqnarray}
\

The case $A_1$ is exceptional. 
For each $N\ge 0$ and 
$\vartheta^{[A_1]}\in X^{[A_1]}$
we define
\begin{eqnarray}
{\bf G}_N^{[{A_1}]}(\vartheta^{[A_1]})
&=&
\prod_{i=0,1} \prod_{a=1}^{N} \int_{|z_{i,a}|=1}\frac{dz_{i,a}}{2\pi \sqrt{-1}z_{i,a}}
\overset{\longrightarrow}{\prod_{~1\leq a \leq N}}
S_0^{[A_1]}(z_{0,a})
\overset{\longrightarrow}{\prod_{~1\leq a \leq N}}
S_1^{[A_1]}(z_{1,a})
\label{def:non-local2}
\\
&\times&
\frac{\displaystyle \prod_{i=0,1} \prod_{1\leq a<b \leq N }
\theta(u_{i,a}-u_{i,b}+\beta)
\theta(u_{i,a}-u_{i,b})}{\displaystyle
\prod_{a=1}^{N}
\prod_{b=1}^{N}
\theta\left(u_{0,a}-u_{1,b}+\beta-\gamma\right)
\theta\left(u_{0,a}-u_{1,b}+\gamma\right)}
\times
\vartheta^{[A_1]}
\left(\sum_{a=1}^{N}u_{0,a},
\sum_{a=1}^{N}u_{1,a}
\right).\nonumber
\end{eqnarray}

\subsection{Commutativity}

\begin{theorem}
\cite{Feigin-Kojima-Shiraishi-Watanabe, Kojima-Shiraishi, Feigin-Jimbo-Mukhin1, Jimbo-Kojima}~~~
In the case of 
${\mathfrak g}=A_l$, $l \geq 1$, and $D_l$, $l \geq 4$, 
the non-local integrals of motion
given in (\ref{def:non-local1}) and (\ref{def:non-local2}) are mutually commutative:
\begin{eqnarray}
[{\bf G}_M^{[{\mathfrak g}]}(\vartheta^{[\mathfrak{g}]}_1), 
{\bf G}_N^{[{\mathfrak g}]}(\vartheta^{[\mathfrak{g}]}_2)]=0
\qquad
\text{for $M,N\ge0$ and 
$\vartheta^{[\mathfrak{g}]}_1, 
\vartheta^{[\mathfrak{g}]}_2\in X^{[\mathfrak{g}]}$}.
\nonumber
\end{eqnarray}
\end{theorem}
~\\
{\bf Conjecture 1.}~
\textit{ 
The same statement holds for ${\mathfrak g}=E_l$, $l=6,7,8$.}

~\\
A word is in order concerning the literature.
In \cite{Feigin-Kojima-Shiraishi-Watanabe} and \cite{Kojima-Shiraishi},
non-local integrals of motion were given for the deformed $W$-algebras of types $A_1$ and $A_l$ $(l \geq 2)$, respectively,
and their commutativity was studied. 
However, to be precise, there were mistakes in the
proof of commutativity of type $A_l$
in \cite{Feigin-Kojima-Shiraishi-Watanabe, Kojima-Shiraishi}.
A corrected proof was  given in \cite{Feigin-Jimbo-Mukhin1}.
Proof of commutativity for type $D_l$ will be given in
\cite{Jimbo-Kojima}.

\section{Discussions}
\label{sec:5}

We end with comments on technical points.

Commutativity of the non-local integrals of motion 
is shown by reducing them to
theta function identities. 
In the case of $A_l$ and $D_l$, 
these theta function identities were shown by mathematical induction 
using the standard
method of complex analysis focusing on poles and zeros.
However, the number of poles and zeros is insufficient to show 
the theta identity corresponding to 
$E_{6,7,8}$ in an analogous way.

In the case of $A_l$,
recently it has become clear that the non-local integrals of motion
can be constructed by trace of the universal $R$ matrix
of the quantum toroidal algebra
\cite{Feigin-Jimbo-Mukhin2}.
According to this construction, commutativity of the non-local integrals of motion follows directly 
from the Yang-Baxter equation.
However, we are not aware of a literature
where the existence of the universal $R$ matrix is established for quantum 
toroidal algebras of type other than $A$. For that reason, 
we chose a direct computational proof of the commutativity.
We remark that a simple minded analog of 
formula (1) for types $B_l$ and $C_l$ does not work,
so some new idea is needed to treat these cases.

\begin{acknowledgement}
TK is thankful for the kind hospitality by the organizing committee of the 16-th International Workshop 
"Lie theory and its application in physics". 
This work is supported by the Grant-in-Aid for Scientific Research
C(19K0350900) and C(25K07041) from Japan Society for Promotion of Science.
\end{acknowledgement}


\begin{thebibliography}{99.}
\bibitem{Gardner} C.S. Gardner, J. Math. Phys. \textbf{12} (1971) 1548-1551
\bibitem{Faddeev-Zakharov} L.D. Faddeev and V.E. Zakharov, Funct. Anal. Appl. \textbf{5} (1971) 280-287
\bibitem{Lax1} P.D. Lax, Commun. Pure and Applied Math. \textbf{21} (1968) 467-490
\bibitem{Lax2} P.D. Lax, Commun. Pure and Applied Math. \textbf{28} (1975) 141-188
\bibitem{Date-Tanaka} E. Date and S. Tanaka,  Prog. Theor. Phys. Suppl. \textbf{59} (1976) 107-125
 \bibitem{Marchenko} V.A. Marchenko, Math. USSR Sbornik \textbf{24} (1974) 319-344
\bibitem{Faddeev-Takhtajan} L.D. Faddeev and L.A. Takhtajan, \textit{Hamiltonian methods in the theory of solitons} (Berlin: Springer-Velag 1987)
\bibitem{Drinfeld-Sokolov} V.G. Drinfel'd and V. V. Sokolov, J. Soviet Math. \textbf{30} (1985) 1975-2036
\bibitem{Feigin-Frenkel1} B.L. Feigin and E. Frenkel, Phys. Lett. \textbf{B246} (1990) 75-81
\bibitem{Kac-Roan-Wakimoto} V.G. Kac, Shi-Shyr Roan, and M. Wakimoto, Commun. Math. Phys. \textbf{241} (1990) 307-342
\bibitem{Feigin-Frenkel2} B.L. Feigin and E. Frenkel, \textit{Infinite analysis} Part \textbf{A, B}, Kyoto, 1991,
Adv. Ser. Math. Phys. {\bf 16} (World Sci. Publ., River Edge, NJ, 1992) pp 197-215
\bibitem{Belavin-Polyakov-Zamolodchikov} A.A. Belavin, A.M. Polyakov and A.B. Zamolodchikov, Nucl. Phys. {\bf B241} (1984) 333-380
\bibitem{Bazhanov-Lukyanov-Zamolodchikov1}
V.V. Bazhanov, S.L. Lukyanov, and A.B. Zamolodchikov, Commun. Math. Phys. {\bf 177} (1996) 381-398
\bibitem{Bazhanov-Lukyanov-Zamolodchikov2}
V.V. Bazhanov, S.L. Lukyanov, and A.B. Zamolodchikov, Commun. Math. Phys. {\bf 190} (1997) 247-278
\bibitem{Bazhanov-Lukyanov-Zamolodchikov3}
V.V. Bazhanov, S.L. Lukyanov, and A.B. Zamolodchikov,
Commun. Math. Phys. {\bf 200} (1999) 297-324 
\bibitem{Bazhanov-Lukyanov-Zamolodchikov4}
V.V. Bazhanov, S.L. Lukyanov, and A.B. Zamolodchikov, 
J. Stat. Phys. {\bf 102} (2001) 567-576
\bibitem{Dorey-Tateo} P. Dorey and R. Tateo, 
J. Phys. {\bf A32} (1999) L419-L425
\bibitem{Bazhanov-Hibberd-Khoroshkin}
V.V. Bazhanov, A.N. Hibberd, and S.M. Khoroshkin, Nucl. Phys. {\bf B622} (2002) 475-547
\bibitem{Feigin-Frenkel3}B.L. Feign and E. Frenkel, {\it Lecture Notes in Math.} {\bf 1920}  (Springer, Berlin, 1996) pp 349-418
\bibitem{Shiraishi-Kubo-Awata-Odake}
J. Shiraishi, H. Kubo, H. Awata, and S. Odake, Lett. Math. Phys. {\bf 38} (1996) 35-51
\bibitem{Awata-Kubo-Odake-Shiraishi} H. Awata, H. Kubo, S. Odake, and J. Shiraishi, Commun. Math. Phys. {\bf 179} (1996) 401-416
\bibitem{Feigin-Frenkel4} B. Feigin and E. Frenkel, Commun. Math. Phys. {\bf 178} (1996) 653-678
\bibitem{Sevostyanov} A. Sevostyanov, Selecta Math. {\bf 8} (2002) 637-703
\bibitem{Frenkel-Reshetikhin} E. Frenkel and N. Reshetikhin, Commun. Math. Phys. {\bf 197} (1998) 1-31
\bibitem{Feigin-Jimbo-Mukhin-Vilkoviskiy} B. Feigin, M. Jimbo,  E.  Mukhin, and I. Vilkoviskiy, 
Selecta Math. (N.S.) {\bf 27} (2021), no. 4, Paper No. 52, 62 pp



\bibitem{Feigin-Kojima-Shiraishi-Watanabe} B. Feigin, T. Kojima, J. Shiraishi, and H. Watanabe, 
\textit{The integrals of motion for the deformed Virasoro algebra}, arXiv:0705.0427v2

\bibitem{Kojima-Shiraishi} T. Kojima and J. Shiraishi, Commun. Math. Phys. {\bf 283} (2008) 795-851

\bibitem{Feigin-Jimbo-Mukhin1} B. Feigin, M. Jimbo, and E. Mukhin, SIGMA {\bf 18} (2022) 051, 31pp

\bibitem{Jimbo-Kojima} M. Jimbo and T. Kojima, Preprint, Now in preparation.

\bibitem{Feigin-Jimbo-Mukhin2}
B. Feigin, M. Jimbo, and E. Mukhin,  J. Phys. {\bf A50}: Math. Theor. (2017) 464001, 28pp 

\end{thebibliography}
\end{document}